\documentclass[reqno,11pt]{amsart}
\newcommand{\be}{\begin{eqnarray}}
\newcommand{\ee}{\end{eqnarray}}
\newcommand{\ben}{\begin{eqnarray*}}
\newcommand{\een}{\end{eqnarray*}}

\newtheorem{thm}{Theorem}[section]

\newtheorem{lema}{Lemma}[section]
\newtheorem{prop}{Proposition}[section]
\usepackage{graphicx}
\usepackage{psfrag}
\usepackage{epsfig,color}

\newcommand{\RR}{\mathbb R}

\begin{document}
\title[Uniqueness of sign changing solutions]
      {On the uniqueness of sign changing bound state solutions of
            a semilinear equation   }\thanks{This research was supported by
        FONDECYT-1070944 for the first author, and
        FONDECYT-1070951  and FONDECYT-1070125 for the second and third author.}
\author[C. Cort\'azar]{Carmen Cort\'azar}
\address{Departamento de Matem\'atica, Pontificia
        Universidad Cat\'olica de Chile,
        Casilla 306, Correo 22,
        Santiago, Chile.}
\email{\tt ccortaza@mat.puc.cl}
\author[M. Garc\'\i a-Huidobro]{Marta Garc\'{\i}a-Huidobro}
\address{Departamento de Matem\'atica, Pontificia
        Universidad Cat\'olica de Chile,
        Casilla 306, Correo 22,
        Santiago, Chile.}
\email{\tt mgarcia@mat.puc.cl}
\author[C. Yarur]{Cecilia S. Yarur}
\address{Departamento de Matem\'atica y C.C.,
        Universidad de Santiago de Chile,
        Casilla 307, Correo 2, Santiago, Chile}
\email{\tt cyarur@usach.cl}



\begin{abstract}
We establish the  uniqueness of the higher  radial bound state solutions
of
$$ \Delta u
+f(u)=0,\quad x\in \RR^n.  \leqno( P) $$ We assume that the nonlinearity  $f\in C(-\infty,\infty)$ is an  odd function satisfying some
convexity and growth conditions, and either has one zero at $b>0$,
is non positive and not identically 0 in $(0,b)$, and  is differentiable and  positive $[b,\infty)$, or is  positive and differentiable in $[0,\infty)$.
\end{abstract}

\maketitle
\today

\section{Introduction and main results}

In this paper we  establish the uniqueness of higher bound state solutions to
$$ \Delta u
+f(u)=0,\quad x\in \RR^n, \leqno( P) $$
in the radial situation. That is, we give conditions on $f$ under which
\begin{eqnarray}\label{eq2}
\begin{gathered}
u''(r)+\frac{n-1}{r}u'(r)+f(u)=0,\quad r>0,\quad n\ge2,\\
u'(0)=0,\quad \lim\limits_{r\to\infty}u(r)=0,
\end{gathered}
\end{eqnarray}
has exactly two solutions, one with $u(0)>0$ and one with with $u(0)<0$,
having a prescribed number of zeros.

Any nonconstant solution to \eqref{eq2} is called a bound state solution. Bound state solutions such that $u(r)>0$ for all $r>0$, are referred to as a first bound state solution, or  a ground state solution. The uniqueness of the first bound state solution of \eqref{eq2} or for the quasilinear situation involving the $m$-Laplacian operator $\nabla\cdot (|\nabla u|^{m-2}\nabla u)$, $m>1$, has been exhaustively studied during the last thirty years, see for example the works \cite{cl},  \cite{coff},
\cite{cfe1}, \cite{cfe2}, \cite{fls}, \cite{Kw}, \cite{m}, \cite{ms},
\cite{pel-ser1}, \cite{pel-ser2},  \cite{pu-ser}, \cite{st}.

We will assume that the function $f:\RR\to\RR$ is continuous, and  that $f$ satisfies  $(f_1)$-$(f_2)$,  where
\begin{enumerate}
\item[$(f_1)$]$f$ is odd,  $f(0)=0$, and there exist $\beta>b>0$ such that $f(s)>0$ for $s>b$,
$f(s)\le 0$, $f(s)\not\equiv 0$ for $s\in[0,b]$,\footnote {The oddness of $f$ is not essential, this assumption can be relaxed to a sign condition: $f(0)=0$, and there exist $b^+>0>b^-$ such that $f(u)>0$ for $u>b^+$, $f(u)<0$ for $u<b^-$, and
$f(u)\le 0$, $f(u)\not\equiv 0$, for $u\in(0,b^+)$ and $f(u)\ge 0$, $f(u)\not\equiv 0$, for $u\in(b^-,0)$} $F(\beta)=0$,  where
$F(s):=\int_0^sf(t)dt$.
\item[$(f_2)$] $f$ is continuous in $[0,\infty)$, continuously differentiable in $(0,\infty)$ and $f'\in L^1(0,1)$.
\end{enumerate}

Our first result  deals with the uniqueness of the $k$-th bound state in space dimension $1<n\le4$:
\begin{thm}\label{main2}
Let $1<n\le4$, $k\in\mathbb N$, and assume that $f$ satisfies  $(f_1)$-$(f_2)$.  If in addition $f$ satisfies
\begin{enumerate}
\item[$(f_4')$] $\displaystyle\Bigl(\frac{F}{f}\Bigr)'(s)\ge \frac{n-2}{2}$ for all $s>\beta$,
\end{enumerate}
then problem \eqref{eq2} has at most one  solution satisfying $u(0)>0$ which has exactly $k-1$ sign changes in $(0,\infty)$.
\end{thm}

Our second result is a strong improvement of the one in \cite{cghy}:

\begin{thm}\label{main1}
Assume that $f$ satisfies  $(f_1)$-$(f_2)$. If $f$ satisfies
\begin{enumerate}
\item[$(f_3)$] $f(s)\ge f'(s)(s-\beta)$, for all $s\ge \beta$, and
\item[$(f_4)$] $\displaystyle\Bigl(\frac{F}{f}\Bigr)'(s)\ge \frac{n-2}{2n}$ for all $s>\beta$,
\end{enumerate}
then problem \eqref{eq2} has at most one  solution satisfying $u(0)>0$ which has exactly one sign change in $(0,\infty)$.
The same conclusion holds if instead of $(f_3)$-$(f_4)$, $f$ satisfies
\begin{enumerate}
\item[$(f_5)$] $\displaystyle \frac{sf'(s)}{f(s)}$ decreases for all $s\ge \beta$, and
\item[$(f_6)$] $\displaystyle \frac{\beta f'(\beta)}{f(\beta)}\le\frac{n}{n-2}$, with $n>2$.
\end{enumerate}

\end{thm}

This work can be seen as a natural continuation of \cite{cghy}, where we established uniqueness of the second bound state solution in the superlinear case.

To the best of our knowledge, there is only one  work (besides \cite{cghy}) concerning the uniqueness of higher  bound states: Troy, see \cite[Theorems 1.1, Theorem 1.3]{troy} studied the existence and uniqueness of the solution to
\eqref{eq2} having exactly one sign change in dimension $n=3$ for
$$f(s)=\begin{cases}
s+1,\quad s\le -1/2,\\ -s,\quad s\in(-1/2,1/2),\\ s-1,\quad s\ge 1/2.
\end{cases}$$
Note that in this case $b=1$, $\beta=1+\sqrt{2}/2$, and for $s>\beta$,
$$(s-\beta)f'(s)=s-1-\frac{\sqrt{2}}{2}<s-1.$$
Hence all assumptions $(f_1)$-$(f_3)$ are satisfied.
Moreover, also $(f_4')$ is satisfied, since for $s\ge\beta$,
$$\Bigl(\frac{F}{f}\Bigr)'(s)=\frac{1}{2}+\frac{1}{4(s-1)^2}\ge\frac{1}{2}=\frac{n-2}{2}\Bigm|_{n=3}\quad\mbox{for all $s\ge\beta$.}$$
Hence, according to our Theorem \ref{main2}, in this case problem \eqref{eq2} has at most one solution with exactly $k$ zeros in $(0,\infty)$ for any $k\in\mathbb N$. Other typical example of a function $f$ satisfying the assumptions of Theorem \ref{main2} is
$$f(s)=s^p-s^q,\quad p>q>0,$$
with no other restriction if $n=2$, and $p^2+q^2\le 1$ when $n=3$.
\medskip

We also deal with the Dirichlet problem in a given ball. In this case we establish non uniqueness of solutions for some $f$ satisfying $(f_1)$-$(f_3)$ (see section \ref{b0}) and we are led to study the situation in the case that $b=0$, that is, $f$ is positive in $(0,\infty)$. More precisely, we assume
\begin{enumerate}
\item[$(f_1')$]$f(0)=0$, and  $sf(s)>0$ for $s>0$,
\item[$(f_2')$] $f$ is continuous in $[0,\infty)$,  continuously differentiable in $(0,\infty)$ and $f'\in L^1(0,1)$,
\item[$(f_3')$] $f(s)\ge sf'(s)$, for all $s>0$,  and 
for any $\varepsilon>0$ there exists $s\in(0,\varepsilon)$ such that $f(s)>sf'(s)$.
\end{enumerate}

\noindent We have imposed the second part in $(f_3')$  to avoid  $f$ linear, for in this case we obviously do not have uniqueness.

On the other hand, it can be shown, see section \ref{b0},
that under these assumptions there do not exist nontrivial bound states,
hence for a given $\rho>0$, we study the Dirichlet problem
\begin{eqnarray}\label{eq2d}
\begin{gathered}
u''(r)+\frac{n-1}{r}u'(r)+f(u)=0,\quad r\in(0,\rho),\quad n\ge2,\\
u'(0)=0,\quad u(\rho)=0,
\end{gathered}
\end{eqnarray}
and prove the following result:
\begin{thm}\label{main2d0}
Assume that $f$ satisfies  $(f_1')$-$(f_3')$, and let $k\in\mathbb N$. Then problem \eqref{eq2d} has at most one  solution satisfying $u(0)>0$ which has exactly $k$ zeros in $(0,\rho)$.
\end{thm}

The existence of sign changing bound state solutions of \eqref{eq2} has been established by Coffman in \cite{coff2} and Mc Leod, Troy and Weissler in \cite{mtw}, where $f:\mathbb R\to\mathbb R$ is locally Lipschitz continuous and satisfies appropriate sign
conditions and is of subcritical growth.  Their proof uses
shooting techniques and a scaling argument. Here we also establish existence by adapting some results in \cite{fls}.
 In \cite{mtw} the function $f$ is assumed to satisfy (besides $(f_1)$ and $(f_2)$)
$$f(u)=C |u|^{p-1}u+g(u),\quad u>0,$$
where $C$ is a positive constant, $g(u)=o(u^p)$ as $u\to\infty$, and $1<p<\displaystyle\frac{n+2}{n-2}$, i.e., it is superlinear and subcritical. They also establish existence for the Dirichlet problem in a ball.

Finally we describe our approach.
In order to prove our results, and due to the oddness of $f$, we will study the behavior of the solutions to the initial value problem
\begin{equation}\label{ivp}
\begin{gathered}
u''(r)+\frac{n-1}{r}u'(r)+f(u)=0
\quad r>0,\quad n\ge2,\\
u(0)=\alpha\quad u'(0)=0
\end{gathered}\end{equation}
for $\alpha\in(0,\infty)$. As usual, we will denote by $u(r,\alpha)$ a $C^2$ solution of \eqref{ivp}.

Our theorems will follow after a series of comparison results between two solutions to \eqref{ivp} with initial value in some small neighborhood of $\alpha^*$, where $u(\cdot,\alpha^*)$ is  a $k$-th bound state, that is, $u(r,\alpha^*)$ is a solution to \eqref{ivp} which has exactly $k-1$ sign changes in $(0,\infty)$ and $\lim\limits_{r\to\infty}u(r,\alpha^*)=0$. We will  show, (see Proposition \ref{delta0} and Lemma \ref{delta1k}), that there exists a neighborhood $V$ of $\alpha^*$ such that any solution to \eqref{ivp} with $\alpha\in V$ has $k$ extremal points in some closed interval $[0,A]$, $A>0$, having extremal values $|E|>\beta$.   In Section \ref{balpha} we follow the ideas of Coffman, see \cite{coff},  and use the function $\varphi(r,\alpha)=\frac{\partial}{\partial\alpha}u(r,\alpha)$ to study the behavior of the solutions between two consecutive extremal points.  In Section \ref{weak} we prove Theorem \ref{main2} through a careful analysis of the behavior of two solutions $u_1(r)=u(r,\alpha_1)$, $u_2(r)=u(r,\alpha_2)$ for $\alpha_1,\ \alpha_2$ in a small neighborhood of $\alpha^*$.  The main tool we use is
the functional
$$Q(s,\alpha)=-4\frac{F}{f}(s)\frac{r(s,\alpha)}{r'(s,\alpha)}-\frac{r^2(s,\alpha)}{(r'(s,\alpha))^2}
-2r^2(s,\alpha)F(s)+H(s),\quad s\neq b,$$
where $H(s)$ is chosen appropriately so that
$$Q'(s,\alpha)=\frac{\partial Q}{\partial s}(s,\alpha)=\Bigl(2(n-2)-4\Bigl(\frac{F}{f}\Bigr)'(s)\Bigr)\frac{r(s,\alpha)}{r'(s,\alpha)},
$$
and the functional $W$ defined by
$$W(s,\alpha)= r(s,\alpha)\sqrt{(u'( r(s,\alpha),\alpha))^2+2F(s)},\quad s\in[U_m(\alpha), \alpha],$$
introduced  in \cite{fls}. Here $ r(s,\alpha)$ denotes the inverse of $u$ between two consecutive extremal points.  In view of hypothesis $(f_4')$, the functional $Q$ allows us to prove some key comparison results concerning the solutions $u_1$ and $u_2$ between their $i-1$-th and $i$-th extremal points, for any $i=1,...,k-1$.

Section \ref{strong} is devoted to the proof of Theorem \ref{main1}, where we use ideas of Pucci, Serrin and Tang in \cite{pu-ser, st} to study the behavior of the solutions in the interval $[\bar U_1,-\beta]$ before the minimum. We do so by considering the celebrated functional introduced first by Erbe and Tang in \cite{et}:
$$P(s,\alpha)=-2n\frac{F}{f}(s)\frac{r^{n-1}(s,\alpha)}{r'(s,\alpha)}-\frac{r^n(s,\alpha)}{(r'(s,\alpha))^2}
-2r^n(s,\alpha)F(s),\quad s\neq b,$$
and the modified functional $\tilde W$ defined by
$$\tilde W(s,\alpha)=r^{n-1}(s,\alpha)\sqrt{(u'(r(s,\alpha),\alpha))^2+2F(s)},\quad s\in[\overline{U}(\alpha), \alpha]$$
where $r(s,\alpha)$ denotes the inverse of $u$ before the first minimum point.

Finally in section 5 we treat the Dirichlet problem and sketch the proof of Theorem \ref{main2d0}.

\section{Preliminaries}\label{prel}

The aim of this section is to establish several properties of the solutions to the initial value problem \eqref{ivp}.

The
functional
\begin{eqnarray}\label{funct-1}
I(r,\alpha)=(u'(r,\alpha))^2+2F(u(r,\alpha)) \end{eqnarray}
will play a fundamental role. A
simple calculation yields \begin{eqnarray}\label{I'}
I'(r,\alpha)=-\frac{2(n-1)}{r}(u'(r))^2,
\end{eqnarray}
and therefore, as $n\ge2$, we have that $I$ is decreasing in $r$.
It can be seen that for $\alpha\in(b,\infty)$, one has $u(r,\alpha)>0$ and $u'(r,\alpha)<0$ for
$r$ small enough, and thus  we can define the extended real number
$$Z_1(\alpha):=\sup\{r>0\ |\ u(s,\alpha)>0\mbox{ and }u'(s,\alpha)<0\ \mbox{ for all }s\in(0,r)\}.$$

Following \cite{pel-ser1}, \cite{pel-ser2} we set
\begin{eqnarray*}
{\mathcal N_1}&=&\{\alpha>0\ :\ u(Z_1(\alpha),\alpha)=0\quad\mbox{and}\quad u'(Z_1(\alpha),\alpha)<0\}\\
{\mathcal G_1}&=&\{\alpha>0\ :\ u(Z_1(\alpha),\alpha)=0\quad\mbox{and}\quad u'(Z_1(\alpha),\alpha)=0\}\\
{\mathcal P_1}&=&\{\alpha>0\ :\ u(Z_1(\alpha),\alpha)>0\}.
\end{eqnarray*}

As in \cite{cfe1}, the sets ${\mathcal N_1}$ and ${\mathcal P_1}$ are open intervals, and moreover, if $\mathcal N_1\not=\emptyset$, then $\mathcal N_1=(a,\infty)$ for some $a>0$. If our problems  have a solution, then ${\mathcal N_1}\not=\emptyset$.
Let
$$\widetilde{\mathcal F}_2=\{\alpha\in\mathcal N_1\ :\ u'(r,\alpha)<0\quad\mbox{for all }r>Z_1(\alpha)\}.$$
For $\alpha\not\in\widetilde{\mathcal F}_2$ we define
$$T_1(\alpha):=\inf\{r>Z_1(\alpha)\ :\ u'(r,\alpha)=0\},\quad \overline{U}_1(\alpha)=u(T_1(\alpha),\alpha),$$
and if $\alpha\in\widetilde{\mathcal F}_2,$ we set $T_1(\alpha)=\infty$. Also, for $\alpha \in \mathcal N_1\setminus \widetilde{\mathcal F}_2$ we can define the extended real number
$$Z_2(\alpha):=\sup\{r>T_1(\alpha)\ |\ u(s,\alpha)<0\mbox{ and }u'(s,\alpha)>0\ \mbox{ for all }s\in(T_1(\alpha),r)\},$$
and set $ U_2(\alpha):=u(Z_2(\alpha),\alpha)=\lim\limits_{r\uparrow Z_2(\alpha)}u(r,\alpha)$.

Let now
$${\mathcal F_2}=\{\alpha\in\mathcal N_1\setminus \widetilde{\mathcal F}_2\ :\ u(Z_2(\alpha),\alpha)<0\},$$
\begin{eqnarray*}
{\mathcal N_2}&=&\{\alpha\in\mathcal N_1\setminus \widetilde{\mathcal F}_2\ :\ u(Z_2(\alpha),\alpha)=0\quad\mbox{and}\quad u'(Z_2(\alpha),\alpha)>0\},\\
{\mathcal G_2}&=&\{\alpha\in\mathcal N_1\setminus \widetilde{\mathcal F}_2\ :\ u(Z_2(\alpha),\alpha)=0\quad\mbox{and}\quad u'(Z_2(\alpha),\alpha)=0\},\\
{\mathcal P_2}&=&\widetilde{\mathcal F}_2\cup
{\mathcal F_2}.
\end{eqnarray*}

For $k\ge 3$, and if ${\mathcal N_{k-1}}\not=\emptyset$, we set
$$\widetilde{\mathcal F}_k=\{\alpha\in\mathcal N_{k-1}\ :\ (-1)^ku'(r,\alpha)<0\quad\mbox{for all }r>Z_{k-1}(\alpha)\}.$$
For $\alpha\not\in \widetilde{\mathcal F}_k$, we set
$$T_{k-1}(\alpha):=\inf\{r>Z_{k-1}(\alpha)\ :\ u'(r,\alpha)=0\},\quad \overline{U}_{k-1}(\alpha)=u(T_{k-1}(\alpha),\alpha),$$
and if $\alpha\in \widetilde{\mathcal F}_k$, we set $T_{k-1}(\alpha)=\infty$. Next, for $\alpha\in \mathcal N_{k-1}\setminus \widetilde{\mathcal F}_k$, we define the extended real number
\begin{eqnarray*}
Z_k(\alpha):=\sup\{r>T_{k-1}(\alpha)\ |\ (-1)^ku(s,\alpha)<0\mbox{ and }(-1)^ku'(s,\alpha)>0\ \\
\mbox{ for all }s\in(T_{k-1}(\alpha),r)\},
\end{eqnarray*}
we set $ U_k(\alpha):=u(Z_k(\alpha),\alpha)=\lim\limits_{r\uparrow Z_k(\alpha)}u(r,\alpha)$. Finally we set
$${{\mathcal F}_k}=\{\alpha\in\mathcal N_{k-1}\setminus \widetilde{\mathcal F}_k\ :\ (-1)^ku(Z_k(\alpha),\alpha)<0\},$$
\begin{eqnarray*}
{\mathcal N_k}&=&\{\alpha\in\mathcal N_{k-1}\setminus \widetilde{\mathcal F}_k\ :\ u(Z_k(\alpha),\alpha)=0\quad\mbox{and}\quad (-1)^ku'(Z_k(\alpha),\alpha)>0\},\\
{\mathcal G_k}&=&\{\alpha\in\mathcal N_{k-1}\setminus \widetilde{\mathcal F}_k\ :\ u(Z_k(\alpha),\alpha)=0\quad\mbox{and}\quad u'(Z_k(\alpha),\alpha)=0\},\\
{\mathcal P_k}&=&\widetilde{\mathcal F}_k\cup
{{\mathcal F}_k}.
\end{eqnarray*}
Concerning the sets ${\mathcal N_k}$ and ${\mathcal P_k}$ we have:

\begin{prop}
\label{open-sets}
The sets ${\mathcal N_k}$ and ${\mathcal P_k}$ are open.
\end{prop}
\begin{proof}
The proof that ${\mathcal N_k}$ is open is by continuity and follows as in \cite{cfe2} with obvious modifications, so we omit it.
\medskip

The proof that $\mathcal P_k$ is open is based in the fact that the functional $I$ defined in \eqref{funct-1} is
decreasing in $r$, and $\alpha\in\mathcal P_k$ if and only if $\alpha\in \mathcal N_{k-1}$ and $I(r_1,\alpha)<0$ for some $r_1\in(0,T_{k-1}(\alpha))$.

 Let $\alpha\in \mathcal P_k$ and assume first that $Z_{k}(\alpha)=\infty$.
 We   claim that $$\lim_{r\to \infty}u(r,\alpha)=-b,\quad \lim_{r\to \infty}u'(r,\alpha)=0.$$
 Since  $u(\cdot,\alpha)$ is monotone ( for all $r>Z_{k-1}$ if   $\alpha\in\widetilde{\mathcal F}_k$ or in
 $(T_{k-1}(\alpha),\infty)$ if $\alpha\in{\mathcal F}_k$), there exists $L$ such that $\lim_{r\to\infty}u(r,\alpha)=L$.
 Furthermore, since  $I(\cdot,\alpha)$ is  decreasing and bounded and $F(s)\to\infty$ as $s\to\pm\infty$, we have that $L$ is finite and
$ \lim_{r\to\infty}u'(r,\alpha)=0.$
   Moreover, from the equation and applying L'H\^opital's rule twice, we conclude that
 $$
 0=\lim_{r\to\infty}\frac{u(r,\alpha)-L}{r^2}= \lim_{r\to\infty}\frac{r^{n-1}u'(r,\alpha)}{2r^n}=-\frac{f(L)}{2n},
 $$
 Thus, $L=-b$ as we claimed, implying that
 $$\lim_{r\to T_{k-1}(\alpha)}I(r,\alpha)=2F(-b)<0.$$
Assume next $Z_k(\alpha)<\infty$ and hence  $\alpha\in{\mathcal F}_k$. Then $T_{k-1}(\alpha)$ is a either a maximum point or a minimum point of $u(\cdot,\alpha)$ implying that either $$0\le -u''(T_{k-1}(\alpha),\alpha)=f(u(T_{k-1}(\alpha),\alpha))$$ and thus $-b<u(T_{k-1}(\alpha),\alpha)<0$ or $$0\ge -u''(T_{k-1}(\alpha),\alpha)=f(u(T_{k-1}(\alpha),\alpha))$$ and thus $0<u(T_{k-1}(\alpha),\alpha)<b$, ($u(T_{k-1}(\alpha),\alpha)\not= \pm b$ from the uniqueness of the  solutions and since $u(0,\alpha)=\alpha)$).
 Hence
 $$I(T_{k-1}(\alpha),\alpha)=2F(u(T_{k-1}(\alpha),\alpha))<0.$$

Conversely, if $\alpha\not\in{\mathcal P_k}$ and $\alpha\in{\mathcal N_{k-1}}$, then $\alpha\in{\mathcal G_k}\cup{\mathcal N_k}$, and thus the
claim follows from the fact that $I(r,\alpha)\ge I(Z_{k}(\alpha),\alpha)\ge 0$ for all $r\in (0, Z_{k}(\alpha)).$
Hence the openness of ${\mathcal P_k}$ follows
from the
continuous dependence of solutions to \eqref{ivp} in the initial value $\alpha$ and from the openness of ${\mathcal N_{k-1}}$.
\end{proof}

Finally in this section we establish the existence of a neighborhood of $\alpha^*$ so that solutions with initial value in this interval cannot be decreasing for all $r>0$.
\begin{prop}\label{delta0}Let $\alpha^*\in\mathcal G_k$, $k\ge 2$. Then there exists  $\delta_0>0$ such that $(\alpha^*-\delta_0, \alpha^*+\delta_0)\subseteq \mathcal{N}_{k-1}\setminus \widetilde{\mathcal F}_k$.
\end{prop}
\begin{proof}Since $\alpha^*\in\mathcal G_k,$
there exists $\tau>T_{k-1}(\alpha^*)$ such that $(-1)^{k}u'(\tau,\alpha^*)>0$. By continuity, there exists $\delta_0>0$ such that
$$(-1)^{k}u'(\tau,\alpha)>0\quad\mbox{for all }\alpha\in(\alpha^*-\delta_0,\alpha^*+\delta_0),$$
implying that
$$T_{k-1}(\alpha)<\tau\quad\mbox{for all }\alpha\in(\alpha^*-\delta_0,\alpha^*+\delta_0),$$ and thus $$(\alpha^*-\delta_0,\alpha^*+\delta_0)
\subset \mathcal N_{k-1}\setminus \widetilde{\mathcal F}_k.$$
\end{proof}

\section{Behavior of the function $\varphi(r,\alpha)=\frac{\partial}{\partial\alpha}u(r,\alpha)$}\label{balpha}
 We will study the behavior of the solutions to the initial value problem \eqref{ivp}. To this end, $\alpha^*\in\mathcal G_k$ is fixed and  $\alpha\in(\alpha^*-\delta_0,\alpha^*+\delta_0)$, where $\delta_0>0$ is given in Proposition \ref{delta0}.

  Under assumptions
 $(f_1)$ and $(f_2)$, the
 functions $u(r,\alpha)$ and $u'(r,\alpha)=\frac{\partial u}{\partial r}(r,\alpha)$ are of class
 $C^1$ in $(0,\infty)\times(b,\infty)$.
We set
$$\varphi(r,\alpha)=\frac{\partial u}{\partial\alpha}(r,\alpha),\quad
'=\frac{\partial }{\partial r}.$$
Then, for any $r>0$ such that $u(r)\neq0$, $\varphi$ satisfies the linear
differential equation
\begin{eqnarray}\label{varphi-eq}
\begin{gathered}
\varphi''(r)+\frac{n-1}{r}\varphi'(r)+ f'(u)\varphi=0,
\quad n\ge2,\\
\varphi(0,\alpha)=1\quad \varphi'(0,\alpha)=0.
\end{gathered}\end{eqnarray}

Set
$$u(r)=u(r,\alpha),\qquad\varphi(r)=\varphi(r,\alpha).$$

\begin{prop}\label{varphi2}
Let $f$ satisfy $(f_1)$-$(f_2)$. Then (i) between two consecutive zeros $r_1< r_2$ of $u'$ there is at least one zero $r^*\in(r_1,r_2)$ of $\varphi$. (ii) Furthermore, if  $\alpha\in\mathcal G_k$, then   $\varphi$ has at least one  zero in $(T_{k-1}(\alpha),Z_k(\alpha))$.
\end{prop}
\begin{proof}
Let $r_1<r_2$ be two consecutive finite zeros of $u'$ (hence $u$ has at most one zero in $(r_1,r_2)$) and assume by contradiction that $\varphi(r)$ does not change sign in $(r_1,r_2)$. Since $u\in C^2(0,\infty)$ and $\varphi\in C^1(0,\infty)$, by differentiating the equation in \eqref{eq2} we obtain that $v=u'$ and $\varphi$ satisfy
\begin{equation}\label{v}
v''+\frac{n-1}{r}v'+\Bigl(f'(u)-\frac{n-1}{r^2}\Bigr)v=0,
\end{equation}
and
\begin{equation}\label{var}\varphi''+\frac{n-1}{r}\varphi'+f'(u)\varphi=0,
\end{equation}
for all $r$ such that $u(r)\neq0$. Hence multiplying \eqref{v} by $r^{n-1}\varphi$ and \eqref{var} by $r^{n-1}v$ and substracting, we obtain
\begin{equation}\label{resta}
(r^{n-1}(v'\varphi-v\varphi'))'(r)=(n-1)r^{n-3}v\varphi.
\end{equation}
Assume first that $v,\varphi>0$ in $(r_1,r_2)$. Integrating \eqref{resta} over $(r_1,r_2)$ we find that
$$r_2^{n-1}v'(r_2)\varphi(r_2)>r_1^{n-1}v'(r_1)\varphi(r_1),$$
a contradiction with the fact that from our choice of the sign for $v$, it must be that $v'(r_2)<0$ and $v'(r_1)>0$. (If $u(\bar r)=0$ for some $\bar r\in(r_1,r_2)$, we integrate \eqref{resta}  over $(r_1,\bar r-\varepsilon)$ and over $(\bar r+\varepsilon, r_2)$, use the continuity of $v$, $v'$, $\varphi$ and $\varphi'$, and then let $\varepsilon\to0$ to obtain a contradiction). Hence $\varphi$ must have a first zero in $(r_1,r_2)$. If either $v$ or $\varphi$ are negative in $(r_1,r_2)$ the proof follows with obvious modifications.

Let now  $\alpha\in\mathcal G_k$. If $Z_k(\alpha)<\infty$, the claim follows from $(i)$. If $Z_k(\alpha)=\infty$,  assume by contradiction that $\varphi$ does not change sign in $(T_{k-1}(\alpha),\infty)$. We may assume without loss of generality that $u'(r)>0$  and $\varphi(r)>0$ for  all $r\in(T_{k-1}(\alpha),\infty)$. From $u'(r)>0$ for  all $r\in(T_{k-1}(\alpha),\infty)$, and $u(r)\to0$ as $r\to \infty$, we find that there exists $r_0>T_{k-1}(\alpha)$ such that $-b<u(r)<0$ for all $r\in(r_0,\infty)$ implying
$$(r^{n-1}u')'=-r^{n-1}f(u)\le0.$$
Thus $r^{n-1}u'$ decreases in $(r_0,\infty)$ implying that
\begin{equation}\label{limitv}
\lim\limits_{r\to \infty}r^{n-1}u'(r)=L\in[0,\infty).
\end{equation} From the equation we find that
$$u''(r)=-\frac{n-1}{r}u'(r)-f(u(r))<0\quad\mbox{for all $r\in (r_0,\infty)$},$$
and thus $v'=u''<0$ for all $r\in (r_0,\infty)$. On the other hand, integrating \eqref{resta} over $(T_{k-1}(\alpha),r)$, for $r\in( r_0,\infty)$, we find that
\ben
r^{n-1}(v'\varphi-v\varphi')(r)=(T_{k-1}(\alpha))^{n-1}v'(T_{k-1}(\alpha))\varphi(T_{k-1}(\alpha))\\
+(n-1)\int_{T_{k-1}(\alpha)}^rt^{n-3}v(t)\varphi(t)dt\\
\ge (n-1)\int_{T_{k-1}(\alpha)}^{r_0}t^{n-3}v(t)\varphi(t)dt=c_0>0\een
 for some positive constant $c_0$. Hence,
$$0>r^{n-1}v'(r)\varphi(r)>r^{n-1}v(r)\varphi'(r)+c_0,$$
which from \eqref{limitv} implies that $\varphi'(r)\le - c_0/(r^{n-1}v)\le -c$ for some positive constant $c$ and therefore
$$\varphi(r)\le \varphi(r_0)-c(r-r_0)\to-\infty\quad\mbox{as }r\to\infty,$$
a contradiction.
\end{proof}

\begin{prop}\label{varphi1}
Let $f$ satisfy $(f_1)$-$(f_3)$. Then $\varphi$ is strictly positive
in $(0,  r(\beta,\alpha))$. 

\end{prop}

\begin{proof}

Multiplying the equation in \eqref{varphi-eq} by $r^{n-1}(u-\beta)$ and integrating
by parts over $(0,r)$, $r\le  r(\beta,\alpha)$, we have that
$$-\int_0^{ r}r^{n-1}u'(r)\varphi'(r)dr+
\int_0^{ r}f'(u(r))\varphi(r)(u(r)-\beta)r^{n-1} dr=0,$$
and a second integration by parts yields
\begin{equation}\label{part1}
 \int_0^{ r}\Bigl(f'(u(t))(u(t)-\beta)-f(u(t))\Bigr)
\varphi(t)t^{n-1} dt=r^{n-1}(u'(r)\varphi(r)-\varphi'(r)(u(r)-\beta)).
\end{equation}
Using now that from $(f_3)$, $f'(u(r))(u(r)-\beta)-f(u(r))\le 0$ for $r\in(0, r(\beta,\alpha))$, we have that if $\varphi(r)=0$ for some $r\in(0,r(\beta,\alpha))$, then $-\varphi'(r)(u(r)-\beta)\le 0$, which is a contradiction since $\varphi'(r)<0$ at such point.
\end{proof}

\medskip
Our next result is an improvement of \cite[Lemma 3.1]{cghy}, where we proved it under an additional superlinear growth assumption on $f$.
\begin{prop}\label{varphi3}
Let $f$ satisfy $(f_1)$-$(f_2)$ and $(f_4)$-$(f_5)$. If the first zero $z>0$ of $\varphi$ occurs in $(0,r(\beta,\alpha)]$, then $\varphi(r)<0$ for $r\in(z,r(b,\alpha))$ and $\varphi'(r(b,\alpha))\le0$.
\end{prop}
\begin{proof}
The proof follows step by step the ideas in \cite{cghy}.
Let the first zero $z>0$ of $\varphi$ occur in $(0,r(\beta,\alpha)]$, set $U_z:=u(z)$  and assume $U_z\ge\beta$. We will show that
$$\frac{U_zf'(U_z)}{f(U_z)}>1.$$
 If not, then by $(f_5)$ we have that
 $$\frac{(s-U_z)f'(s)}{f(s)}<\frac{sf'(s)}{f(s)}\le 1\quad\mbox{for all $s\ge U_z$},$$
 and we can argue as in the proof of Proposition \ref{varphi1} (with $\beta$ replaced by $U_z$) to obtain the contradiction
 \ben
 \int_0^{ z}\Bigl(f'(u(t))(u(t)-U_z)-f(u(t))\Bigr)
\varphi(t)t^{n-1} dt=r^{n-1}(u'(r)\varphi(r)-\varphi'(r)(u(z)-U_z))=0.
\een
We conclude that there exists $c>0$ such that
$$\frac{U_zf'(U_z)}{f(U_z)}=1+\frac{2}{c}.$$
Moreover, from $(f_5)$-$(f_6)$, it must be that $c\ge n-2$.
Then, since by  $(f_5)$, the function
$$r\to  c\frac{u(r)f'(u(r))}{f(u(r))}- c-2$$
is increasing in $(0,  r(b,\alpha))$, we have that
$$\phi(r):= f(u(r))\Bigl( c\frac{u(r)f'(u(r))}{f(u(r))}- c-2\Bigr)$$
is non positive in $(0,z)$ and nonnegative in $(z, r(b,\alpha))$.

Let us set $v(r)=ru'(r)+ cu(r)$. Then $v$ satisfies
$$v''+\frac{n-1}{r}v'+ f'(u(r))v=\phi(r),$$
and, as long as $\varphi(r)$ does not change sign in $(z,r)$, with $r\in(z,r(b,\alpha))$, we have
\begin{eqnarray}\label{41.4}
0&\ge&\int_0^rt^{n-1}\varphi(t)\phi(t)dt
=\int_0^rt^{n-1}(\varphi\Delta v-v\Delta\varphi)dt\nonumber\\
&=&r^{n-1}(\varphi(r)v'(r)-\varphi'(r)v(r)),
\end{eqnarray}
and therefore
\begin{eqnarray}\label{contra2}
\varphi(r)v'(r)-\varphi'(r)v(r)\le 0,
\end{eqnarray}
implying in particular that $v(z)\le 0$.
On the other hand, using that $c\ge n-2$ we have that
$$v'(r)=ru''(r)+(c+1)u'(r)\le  ru''(r)+(n-1)u'(r)=-rf(u(r))<0$$
for all $r\in(0,r(b,\alpha))$.
Now we can prove that $z$ is the only zero of $\varphi$ in $(0,r(b,\alpha))$. Indeed, if $\varphi$ has a second zero at $z_1\in r(b,\alpha))$, then from \eqref{contra2}, it must be that $v(z_1)\ge 0$, contradicting $v'(r)<0$ in $(0,r(b,\alpha))$. Hence $\varphi$ has exactly one zero in $(0,r(b,\alpha)]$.

Finally, evaluating \eqref{contra2} at $r= r(b,\alpha)$,
we find that
$$\varphi( r(b,\alpha))v'( r(b,\alpha))-\varphi'( r(b,\alpha))v( r(b,\alpha))\le 0,$$
implying $\varphi'( r(b,\alpha))\le 0$.
\end{proof}


\section{Uniqueness of bound states}\label{weak}
Assume that  $\alpha^*\in\mathcal G_k$. The following result  deals with
 the existence of  a neighborhood $V$ of $\alpha^*$ such that any solution to \eqref{ivp} with $\alpha\in V$
has its  minimum values  satisfying $U<-\beta$ and its maximum values satisfying $U>\beta.$

We observe that $u(\cdot,\alpha)$ is invertible in each interval $(T_{i-1}(\alpha),T_i(\alpha))$, $T_0(\alpha)=0$, $i=1,2,\ldots,k-1$, and we  denote by  $r(\cdot,\alpha)$ its inverse at the intervals where $u$ decreases and by $\bar r(\cdot,\alpha)$ its inverse at intervals where $u$ increases.

\begin{lema}\label{delta1k} Let $f$ satisfy $(f_1)$-$(f_2)$, and let $\alpha^*\in\mathcal G_k$. Then, there exist $a>0$ and $\delta_1>0$,
such that for any $\alpha\in(\alpha^*-\delta_1,\alpha^*+\delta_1)$, $u(\cdot, \alpha)$ has exactly $k$ extremal points in $[0, T_{k-1}(\alpha^*)+a]$. The extremal values $E$ of $u(\cdot,\alpha)$ satisfy
$E<-\beta$ if $E $ is a minimum value, while $E>\beta$ if $E$  a maximum value. Moreover, if $\alpha_1<\alpha_2$ are two values in
$(\alpha^*-\delta_1,\alpha^*+\delta_1)$, then
\begin{enumerate}
\item[(i)]  the corresponding solutions $u_1$ and $u_2$ intersect between any two of their consecutive extremal points, and
\item[(ii)] there exists an intersection point in $(T_{k-1}(\alpha^*),Z_k(\alpha^*))$.
    \end{enumerate}
\end{lema}
\begin{proof}
Let $\delta_0$ be given as in Proposition \ref{delta0}. The assumption $\alpha^*\in\mathcal G_k$ implies that the functional defined in \eqref{funct-1} satisfies
$$I(Z_k(\alpha^*),\alpha^*)= 0,$$
and thus $I(r,\alpha^*)>0$ for all $r\in(0,Z_k(\alpha^*))$. In particular, for any $i=1,2,\ldots, k-1$, we have
$$2F(u(T_{i}(\alpha^*),\alpha^*))=I(T_{i}(\alpha^*),\alpha^*)>0,$$
implying that $|u(T_{i}(\alpha^*),\alpha^*)|>\beta$.
Hence, from the continuity of $u$ and $T_i(\alpha)$ for $\alpha\in(\alpha^*-\delta_0,\alpha^*+\delta_0)$,
we conclude that there exists $\bar\delta_1<\delta_0$ such that the first assertion of the lemma holds.

From Proposition \ref{varphi2}, for each $i=1,2,\ldots,k-1$, there exists $r^*\in(T_{i-1}(\alpha^*),T_i(\alpha^*))$ such that $\varphi(r^*,\alpha^*)=0$. Hence without loss of generality we may assume that  there exist $r^-<r^*<r^+$ such that $\varphi(r^+,\alpha^*)<0<\varphi(r^-,\alpha^*)$. By continuity, there exists $\delta_1\in(0,\bar\delta_1)$ such that $\varphi(r^-,\alpha)>0$ and $\varphi(r^+,\alpha)<0$ for all $\alpha\in(\alpha^*-\delta_2,\alpha^*+\delta_2)$. Since
$$u(r,\alpha_2)-u(r,\alpha_1)=\int_{\alpha_1}^{\alpha_2}\varphi(r,\alpha)d\alpha,$$
which is positive at $r=r^-$ and negative at $r=r^+$, and thus $(i)$  is proved. $(ii)$ follows in the same way.

\end{proof}

\subsection{Proof of Theorem \ref{main2}}
\mbox{ }\\

We recall that in Theorem \ref{main2}, $2\le n\le 4$.
Let $m <M$ be such that $r(s,\alpha)$ is defined and decreasing in $[m, M]$.
For $s\in[m,M]$ we set

$$Q(s,\alpha)=-4\frac{F}{f}(s)\frac{r(s,\alpha)}{r'(s,\alpha)}-\frac{r^2(s,\alpha)}{(r'(s,\alpha))^2}
-2r^2(s,\alpha)F(s) + H(s),$$ with
$$
H'(s)= -4(n-2)\frac{F}{f}(s).$$
Then,
\begin{equation}\label{qderiv}Q'(s,\alpha)=\frac{\partial Q}{\partial s}(s,\alpha)=\Bigl(2(n-2)-4\Bigl(\frac{F}{f}\Bigr)'(s)\Bigr)\frac{r(s,\alpha)}{r'(s,\alpha)}.
\end{equation}
Similarly, for $\overline{m}<\overline{M}$ such that $\bar r(s,\alpha)$ is defined and increasing in $[\overline{m},\overline{M}]$, we define
$$\bar Q(s,\alpha)=-4\frac{F}{f}(s)\frac{\bar r(s,\alpha)}{\bar r'(s,\alpha)}-\frac{\bar r^2(s,\alpha)}{(\bar r'(s,\alpha))^2}
-2\bar r^2(s,\alpha)F(s) + \bar H(s),$$
with
$$
\bar H'(s)= -4(n-2)\frac{F}{f}(s).$$

Note that if $(f'_4)$  holds, then
$Q'(s,\alpha)\ge 0$ for all $s\in [m, M]$ and $ \bar Q'(s,\alpha)\le 0$
for all $s\in[\overline{m} ,\overline{M}]$.
\medskip

Let now $a$ and $\delta_1$ be as in Lemma \ref{delta1k}, let $\alpha_1,\ \alpha_2\in(\alpha^*-\delta,\alpha^*+\delta)$, with $\alpha_1<\alpha_2$,
and for $j=1,2$ set
$$u_j(r)=u(r,\alpha_j),\quad r_j(s)=r(s,\alpha_j),\quad\mbox{and}\quad
Q_j(s)=Q(s,\alpha_j).$$
 Let
$$M_1,\ m_{1},\mbox{be the $i$-th consecutive local maximum and minumum values of $u_1$}, $$ and
 $$M_2,\ m_{2},\mbox{be the $i$-th  consecutive local maximum and minumum  values of $u_2$} $$
for $r\in[0,T_{k-1}(\alpha^*)+a]$.
The behavior of the solutions for $r>T_{k-1}(\alpha^*)$ will be studied separately.
We have
\begin{prop}\label{qmaxmin}
Assume that $f$ satisfies $(f_1)$-$(f_2)$ and $(f_4')$, and let $\alpha^*\in\mathcal G_k$. Then, there exists $\delta_{2,i}\in(0,\delta_1)$, with $\delta_1$ as in Lemma \ref{delta1k}, such that for any  $\alpha_1,\ \alpha_2\in(\alpha^*-\delta_{2,i},\alpha^*+\delta_{2,i})$ with $\alpha_1<\alpha_2$ we have that if
$$ M_1<M_2\ \mbox{and} \ Q_1(M_1)> Q_2(M_2), $$
then
$$ m_1>m_2\ \mbox{and} \ Q_1(m_1)> Q_2(m_2). $$

\end{prop}

In order to prove this result we need a separation lemma, so for $j=1,2$ we  consider the functional $W_j$ defined below,  introduced  in \cite{ fls}:
$$ W_j(s)=  r_j(s)\sqrt{(u_j'( r_j(s)))^2+2F(s)}, \quad s\in[m_j,M_j],$$
The functional $W_j$ is well defined in this interval, since
$(u_j'(r))^2+2F(u_j(r))>0$ for $r\in[0,T_{k-1}(\alpha^*)+a]$.

\begin{lema}\label{sep4}
Assume that $f$ satisfies $(f_1)$-$(f_2)$, and let $\alpha^*\in\mathcal G_k$. Let $\alpha_1,\ \alpha_2\in(\alpha^*-\delta_1,\alpha^*+\delta_1)$ with $\alpha_1<\alpha_2$ and $\delta_1$ as in Lemma \ref{delta1k}. Assume that there exists $ U\in[-\beta,\beta]$ such that
\begin{equation} \label{barcompa}  r_1(U)\ge  r_2(U)\quad\mbox{and}\quad W_1( U)<W_2(U).\end{equation}
 Then
$$ r_1(s)> r_2(s), \quad W_1(s)<W_2(s),\quad\mbox{ for all } s\in[-\beta, U].$$
\end{lema}

\begin{proof} Clearly, $|r'_1(U)|>|r'_2(U)|$, and thus $r_1>r_2$ in some small left neighborhood of $U.$
Hence, there exists $c\in[-\beta,U)$ such that   $$ W_1\le W_2,  \quad r_1> r_2, \quad\mbox{and} \quad r'_1< r'_2
\quad\mbox{ in $[c,U)$ }.$$

Next, we will show that $ W_1- W_2$ is  increasing in $[c,U)$. This will imply that the infimum of such $c$ is $-\beta$, proving the lemma.

From the definition of $ W_j(s)$ we have
\ben\frac{\partial  W_j}{\partial s}(s)=\frac{-2F(s)+(n-2)(u_j'(r_j(s)))^2}{|u_j'(r_j(s))|\sqrt{(u_j'(r_j(s)))^2+
2F(s)}}.
\een
As $F(s)\le 0$ for $s\in[-\beta,\beta]$, we have that the function
$$h(p)=\frac{-2F(s)}{p\sqrt{p^2+2F(s)}}+\frac{(n-2)p}{\sqrt{p^2+2F(s)}},\quad p>0,$$
is decreasing, and thus, for $s\in [c,U)$, and using that $|u_1'(r_1(s))|<|u_2'(r_2(s))|$, we obtain
$$\Bigl(\frac{\partial  W_1}{\partial s}-\frac{\partial  W_2}{\partial s}\Bigr)(s)=h(|u_1'(r_1(s))|)-h(|u_2'(r_2(s))|)>0$$
as we claimed.

\end{proof}

\begin{proof}[{\bf Proof of Proposition \ref{qmaxmin}}]
First we note that since $Q_2$ is strictly increasing, and $M_1<M_2$, it holds that $Q_1(M_1)>Q_2(M_1)$.

Let $M^*$ denote the $i$-th maximum value of $u(\cdot,\alpha^*)$. Since  $u'( r({M^*},\alpha^*),\alpha^*)=0$ and $4\frac{F}{f}({M^*}) >0$, by continuity there exists  $\delta_{2,i}<\delta_1$ such that for any $\alpha_1,\ \alpha_2\in(\alpha^*-\delta_{2,i},\alpha^*+\delta_{2,i})$, we have
$$4\frac{F}{f}({M}_{1}) >- r_2({M}_{1})u_2'( r_2({M}_{1})),$$
and hence
$$4\frac{F}{f}({M}_{1})r_2({M}_{1})u_2'( r_2({M}_{1}))+( r_2({M}_{1}))^2(u_2'( r_2({M}_{1})))^2<0.$$
Therefore,
\ben
0&<&(Q_1- Q_2)({M}_{1})\\
&=&4\frac{F}{f}({M}_{1}) r_2({M}_{1})u_2'( r_2({M}_{1}))+( r_2({M}_{1}))^2(u_2'( r_2({M}_{1})))^2
+2F({M}_{1})( r_2^2- r_1^2)({M}_{1})\\
&<&2F({M}_{1})( r_2^2- r_1^2)({M}_{1}),
\een
implying
$$ r_1({M}_{1})< r_2({M}_{1}).$$

From Lemma \ref{delta1k} there exists a greatest intersection point $U_I$ of $r_1$ and $r_2$ in $[\max\{m_1,m_2\}, M_1]$.

Let us set $$U=\min\{-\beta, U_I\}.$$ We will show that
 \begin{equation} \label{q0}
(Q_1-Q_2)(U)>0,\qquad \mbox{and} \qquad \frac{r_1}{|r'_1|}(U)< \frac{r_2}{|r'_2|}(U).
\end{equation}

We distinguish the following cases according to the position of $U_I$:

\noindent Case 1.  $U_I\in  [\beta, M_1].$
We will prove first that
$$ \frac{r_1}{|r'_1|}(s)<\frac{r_2}{|r'_2|}(s),\quad \mbox{for all $s\in [U_I, M_1].$} $$
Indeed, since $u'_1(r_1(M_1))=0,$ we have that this inequality holds for $s=M_1$. Assume now that there exists $t\in(U_I,M_1)$ such that
$$ \frac{r_1}{|r'_1|}(s)<\frac{r_2}{|r'_2|}(s),\quad \mbox{for all $s\in (t, M_1)$ and}\quad  \frac{r_1}{|r'_1|}(t)=\frac{r_2}{|r'_2|}(t). $$
As
$$\frac{d}{ds}(\frac{r_1}{|r'_1|}-\frac{r_2}{|r'_2|})(t)=f(t)(r_2|r_2'|-r_1|r_1'|)(t)=f(t)\frac{|r_1'|}{r_1}(t)(r_2^2-r_1^2)(t)>0,
$$
we obtain a contradiction.

Assume next that there exists  $t\in[\beta,U_I)$ such that
$$ \frac{r_1}{|r'_1|}(s)<\frac{r_2}{|r'_2|}(s),\quad \mbox{for all $s\in (t, M_1)$ and}\quad  \frac{r_1}{|r'_1|}(t)=\frac{r_2}{|r'_2|}(t). $$
Then, from $(f_4')$,
$$ (Q_1-Q_2)'(s)=4\Bigl(\frac{r_1}{|r'_1|}(s)-\frac{r_2}{|r'_2|}(s)\Bigr)\Bigl(\Bigl(\frac{F}{f}\Bigr)'(s)-\frac{n-2}{2}\Bigr)<0,\ s\in (t, M_1)
$$
implying that $$ 0>-2F(t)(r_1^2(t)-r_2^2(t))=(Q_1-Q_2)(t)>(Q_1-Q_2)(M_1)>0,$$
a contradiction. We conclude that
$$(Q_1-Q_2)(\beta)>(Q_1-Q_2)(M_1)>0$$
implying
$$ \frac{r_1}{|r'_1|}(\beta)<\frac{r_2}{|r'_2|}(\beta)\quad\mbox{and}\quad r_1(\beta)\ge r_2(\beta).$$
Now we can use Lemma \ref{sep4} with $U=\beta$, to
  obtain that $r_1(-\beta)>r_2(-\beta)$ and $W_1(-\beta)<W_2(-\beta)$, implying \eqref{q0} at $U=-\beta$.
\medskip

\noindent Case 2.  $U_I\in [-\beta, \beta].$ In this case $W_1(U_I)<W_2(U_I)$ and $r_1(U_I)=r_2(U_I)$, hence by Lemma \ref{sep4}, we conclude $W_1(-\beta)<W_2(-\beta)$ implying that \eqref{q0} holds.

\noindent Case 3.  $U_I\in  [\max\{m_1,m_2\}, -\beta].$
In this case it is straightforward to verify that  $$(Q_1-Q_2)(U_I)>0,$$
and hence in this case \eqref{q0} also holds.

To end the proof,  assume that there exists $\tau\in(\max\{m_1,m_2\}, U]$ such that
$$ \frac{r_1}{|r'_1|}(s)<\frac{r_2}{|r'_2|}(s),\quad \mbox{for all $s\in (\tau, U],$} $$ and
$$\frac{r_1}{|r'_1|}(\tau)=\frac{r_2}{|r'_2|}(\tau).$$
Then,
$$ (Q_1-Q_2)'(s)=4\Bigl(\frac{r_1}{|r'_1|}(s)-\frac{r_2}{|r'_2|}(s)\Bigr)\Bigl(\Bigl(\frac{F}{f}\Bigr)'(s)-\frac{n-2}{2}\Bigr)<0,\ s\in (\tau, U]
$$
implying that $$ 0>-2F(\tau)(r_1^2(\tau)-r_2^2(\tau))=(Q_1-Q_2)(\tau)>(Q_1-Q_2)(U)>0,$$
a contradiction, and thus
$$ \frac{r_1}{|r'_1|}(s)<\frac{r_2}{|r'_2|}(s),\quad \mbox{for all $s\in [\max\{m_1,m_2\}, U).$} $$
Therefore,
$$ \max\{m_1,m_2\}=m_1, \quad (Q_1-Q_2)'(s)>0, \quad \mbox{for all} \quad s\in [m_1, U),$$ which yields
$Q_1(m_1)>Q_2(m_1)$. Since $Q_2$  increases and $m_1>m_2$, it follows that $Q_1(m_1)>Q_2(m_2)$, ending the proof of the proposition.
\end{proof}
Similarly we set
$$\bar m_1,\ \bar M_{1}\mbox{ the $i$-th consecutive local   minumum and maximum of $u_1$}, $$ and
 $$\bar m_2,\ \bar M_{2}\mbox{ the $i$-th consecutive local minumum and maximum   of $u_2$}, $$
 for $r\in[0,T_{k-1}(\alpha^*)+a]$.

We have the following result.
\begin{prop}\label{qminmax}
Assume that $f$ satisfies $(f_1)$-$(f_2)$ and $(f_4')$, and let $\alpha^*\in\mathcal G_k$. Then, there exists $\bar\delta_{2,i}\in(0,\delta_1)$, with $\delta_1$ as in Lemma \ref{delta1k}, such that for any  $\alpha_1,\ \alpha_2\in(\alpha^*-\bar\delta_{2,i},\alpha^*+\bar\delta_{2,i})$ with $\alpha_1<\alpha_2$ we have that if
$$ \bar m_1>\bar m_2\ \mbox{and} \ \bar Q_1(\bar m_1)>\bar Q_2(\bar m_2), $$
then
$$ \bar M_1<\bar M_2\ \mbox{and} \ \bar Q_1(\bar M_1)> \bar Q_2(\bar M_2). $$

\end{prop}
\begin{proof}
It follows from Proposition \ref{qmaxmin} considering
$v(r,\alpha_j)=-u(r,\alpha_j)$.
\end{proof}
Combining Propositions \ref{qmaxmin} and \ref{qminmax} we obtain the following result.
\begin{prop}\label{qminmaxmin}
Assume that $f$ satisfies $(f_1)$-$(f_2)$ and $(f_4')$, and let $\alpha^*\in\mathcal G_k$. Let $\delta=\min\limits_{i}\{\delta_{2,i},\bar\delta_{2,i}\}$, and let $\alpha_1,\ \alpha_2\in(\alpha^*-\delta,\alpha^*+\delta)$.
\begin{enumerate}
\item[(i)]
If $k$ is even, then the $k$-th extremal points $T_{k-1}(\alpha_i)$ are minima,
$$m_1>m_2\quad\mbox{and}\quad Q_1(m_1)>Q_2(m_2),$$
where $m_i=u_i(T_{k-1}(\alpha_i))$.
\item[(ii)]If $k$ is odd, then the $k$-th extremal points $T_{k-1}(\alpha_i)$ are maxima,
$$M_1<M_2\quad\mbox{and}\quad Q_1(M_1)>Q_2(M_2),$$
where $M_i=u_i(T_{k-1}(\alpha_i))$.
\end{enumerate}
\end{prop}
\begin{proof}  As $T_0(\alpha_i)=0$ is the first extremal point of $u_i$, we have
$$u_1(T_0(\alpha_1))=\alpha_1<\alpha_2=u_2(T_0(\alpha_2)).$$
Moreover, as $\alpha_i>\beta$, $H$ is decreasing in $[\beta,\infty)$ and therefore
$$Q_1(\alpha_1)=H(\alpha_1)>H(\alpha_2)=Q_2(\alpha_2).$$
Hence, for the first extremal points, the assumption of Proposition \ref{qmaxmin} holds and thus,
$$u_1(T_1(\alpha_1))>u_2(T_1(\alpha_2)),\quad\mbox{and}\quad Q_1(u_1(T_1(\alpha_1)))>Q_2(u_2(T_1(\alpha_2))).$$
Applying alternatively  Proposition \ref{qminmax} and Proposition \ref{qmaxmin} we obtain the result.
\end{proof}
We proceed now to our final step. To this end, we may assume without loss of generality that $k$ is odd, so that $T_{k-1}(\alpha_j)$ is a maximum point, and we fix $\delta$ as given in Proposition \ref{qminmaxmin}.

\begin{prop}\label{sep5q}
Assume that $f$ satisfies $(f_1)$-$(f_2)$ and $(f_4')$, and let $\alpha^*\in\mathcal G_k$. Let $\alpha_1,\ \alpha_2\in(\alpha^*-\delta,\alpha^*+\delta)$ with $\alpha_1<\alpha_2$.

If $\alpha_1\in \mathcal G_k\cup \mathcal N_k$,
then $\alpha_2\in \mathcal N_k$, 
\begin{equation}\label{sep55q}Z_k(\alpha_1)>Z_k(\alpha_2)\quad\mbox{and}\quad |u'_1(Z_k(\alpha_1))|<|u'_2(Z_k(\alpha_2))|.\end{equation}
If $\alpha_2\in\mathcal G_k$, then $\alpha_1\in{\mathcal F_k}$.
\end{prop}
In order to prove this result we need the following separation lemma which can be found in \cite[Lemma 4.4.1]{cghy}. Its proof is very similar to that of Lemma \ref{sep4} and thus we omit it.  Let

  $$S_j:= \inf\{s\in (U_k(\alpha_j), M_j))\ :\ \quad |u_j'( r_j(s))|^2+2F(s)>0\},$$
where $M_j=u_j(T_{k-1}(\alpha_j))$. We note that $S_j=0$ if and only if $\alpha_j\in\mathcal G_k\cup\mathcal N_k$.
\begin{lema}\label{sep4ultimo}
Assume that $f$ satisfies $(f_1)$-$(f_2)$, and let $\alpha^*\in\mathcal G_k$. Let $\alpha_1,\ \alpha_2\in(\alpha^*-\delta,\alpha^*+\delta)$ with $\alpha_1<\alpha_2$.
Assume that there exists $U\in[0,\beta]$ such that
\begin{equation} \label{barcompa2}   r_1(U)\ge   r_2(U)\quad\mbox{and}\quad  W_1( U)< W_2(U).\end{equation}
Then,   $S_1\ge S_2$
 and
$$ r_1(s)> r_2(s), \quad  W_1(s)< W_2(s),\quad\mbox{and \; $|u'_1( r_1(s))|<|u'_2( r_2(s))|$} \quad s\in[S_1,U).$$
\end{lema}
\begin{proof}[{\bf Proof of Proposition \ref{sep5q}.}] Let $r_I$ denote the first intersection point of $u_1$ and $u_2$ in $(T_{k-1}(\alpha^*), Z_k(\alpha^*))$ guaranteed by Lemma \ref{delta1k}(ii) and  $U_I=u_j(r_I)$.
Arguing as in the proof of Proposition \ref{qmaxmin}, cases 1 and 2, this time with $U=\min\{\beta,U_I\}$, we obtain that \eqref{barcompa} holds. Hence, by Lemma \ref{sep4ultimo}, we have
$S_1\ge S_2$,
$$ r_1(s)> r_2(s), \quad  W_1(s)< W_2(s),\quad\mbox{and \; $|u'_1( r_1(s))|<|u'_2( r_2(s))|$} \quad\mbox{for all } s\in[S_1,U).$$
If $\alpha_1\in\mathcal G_k\cup\mathcal N_k$, then $S_1=0$ implying $S_2=0$ and $\alpha_2\in\mathcal G_k\cup\mathcal N_k$. As
$Z_k(\alpha_1)=r_1(0)>r_2(0)=Z_k(\alpha_2)$ and $|u_1'(Z_k(\alpha_1))|<|u_2'(Z_k(\alpha_2))|$ we conclude that $\alpha_2\in\mathcal N_k$.

If $\alpha_2\in\mathcal G_k$, then $S_2=0$. As $|u_2'(Z_k(\alpha_2))|=0$, we conclude that $S_1>0$ implying $\alpha_1\in\mathcal F_k$.
\end{proof}
\begin{proof}[{\bf Proof of Theorem \ref{main2}} ]
Let $\alpha^*\in\mathcal G_k$ hence by Proposition \ref{sep5q}, $(\alpha^*,\alpha^*+\delta)\subset
\mathcal N_k$. Let
$$\bar\alpha=\sup\{\alpha>\alpha^*\ :\ (\alpha^*,\alpha)\subset\mathcal N_k\}.$$
Assume $\bar\alpha<\infty$. Since $\mathcal P_k$ and $\mathcal N_k$ are open, we deduce that $\bar\alpha\in\mathcal G_k$. By Proposition \ref{sep5q}, $(\bar\alpha-\delta,\bar\alpha)\subset{\mathcal F_k}$, a contradiction, and thus $(\alpha^*,\infty)\subset \mathcal N_k$. Hence, there exists at most one solution of \eqref{eq2} with exactly $k-1$ sign changes in $(0,\infty)$.

\end{proof}

\subsection{Proof of Theorem \ref{main1}}\label{strong}

\mbox{ }\\

In what follows we use  the ideas of Pucci, Serrin and Tang in \cite{pu-ser, st}. For
$s\in(\overline{U}_1(\alpha),-\beta]$ we set
$$P(s,\alpha)=-2n\frac{F}{f}(s)\frac{r^{n-1}(s,\alpha)}{r'(s,\alpha)}-\frac{r^n(s,\alpha)}{(r'(s,\alpha))^2}
-2r^n(s,\alpha)F(s).$$
Then,
\begin{equation}\label{pderiv}P'(s,\alpha)=\frac{\partial P}{\partial s}(s,\alpha)=\Bigl(n-2-2n\Bigl(\frac{F}{f}\Bigr)'(s)\Bigr)\frac{r^{n-1}(s,\alpha)}{r'(s,\alpha)}.
\end{equation}

By $(f_4)$ it holds that
$P'(s,\alpha)\ge 0$ for all $s\in(\overline{U}_1(\alpha),-\beta]$.
\medskip

In this case we can prove the analogue of Proposition \ref{qmaxmin} but only for the first maximal and minimal points of $u_1$ and $u_2$.
Let now $\alpha_1,\ \alpha_2\in(\alpha^*-\delta,\alpha^*+\delta)$, with $\alpha_1<\alpha_2$,
and set
$$
P_1(s)=P(s,\alpha_1),\quad P_2(s)=P(s,\alpha_2),$$
$$m_1=u_1(T_1(\alpha_1)),\qquad m_{2}=u_2(T_1(\alpha_2)).$$

We have

\begin{prop}\label{pmaxmin}
Assume that $f$ satisfies $(f_1)$-$(f_3)$ and $(f_4)$, or $(f_1)$-$(f_2)$ and $(f_5)$-$(f_6)$, and let $\alpha^*\in\mathcal G_k$. Let $\alpha_1,\ \alpha_2\in(\alpha^*-\delta,\alpha^*+\delta)$ with $\alpha_1<\alpha_2$ and $\delta=\delta_1$ as in Lemma \ref{delta1k}. Then,
\be\label{ma1} m_1>m_2\ \mbox{and} \ P_1(m_1)> P_2(m_2). \ee

\end{prop}

In order to prove this result we need the following variations of  lemma \ref{sep4}, so for $j=1,2$ we  consider the functional $\tilde W_j$ defined below,
$$\tilde W_j(s)=  r_j^{n-1}(s)\sqrt{(u_j'( r_j(s)))^2+2F(s)}, \quad s\in[m_j,\alpha_j].$$

From Lemma \ref{delta1k}, the solutions $u_1$ and $u_2$ intersect at a first $r_I>0$. Set $U_I= u_1(r_I)=u_2(r_I)$.

\begin{lema}\label{sep1}
Let $f$ satisfy $(f_1)$-$(f_3)$. Let $\alpha_1,\ \alpha_2\in(\alpha^*-\delta,\alpha^*+\delta)$ with $\alpha_1<\alpha_2$ and $\delta=\delta_1$ as in Lemma \ref{delta1k}. If $U_I\in[-\beta,\beta]$
 then
$$r_1( s)> r_2( s)\quad\mbox{and}\quad \tilde W_1( s)<\tilde W_2( s),\quad\mbox{for all }s\in[-\beta,U_I).$$
\end{lema}
\begin{proof} Clearly, $|r'_1(U_I)|>|r'_2(U_I)|$, and thus $r_1>r_2$ in some small left neighborhood of $U_I.$
Hence, there exists $c\in[-\beta,U_I)$ such that   $$\tilde W_1\le\tilde W_2,  \quad r_1> r_2, \quad\mbox{and} \quad r'_1< r'_2
\quad\mbox{ in $[c,U_I)$ }.$$

Next, we will show that $\tilde W_1-\tilde W_2$ is  increasing in $[c,U_I)$. This will imply that the infimum of such $c$ is $-\beta$, proving the theorem.

From the definition of $\tilde W(s, \alpha)$ we have
\ben\frac{\partial \tilde W}{\partial s}(s,\alpha)=\frac{2(n-1)r^{n-2}(s,\alpha)F(s)}{u'(r(s,\alpha),\alpha)\sqrt{(u'(r(s,\alpha),\alpha))^2+
2F(s)}},
\een
and thus, for $s\in [c,U_I)$,
\begin{eqnarray*}
&& \frac{1}{2(n-1)}\Bigl(\frac{\partial \tilde W_1}{\partial s}(s)-\frac{\partial \tilde W_2}{\partial s}(s)\Bigr)\\
&=&F(s)\Bigl(
\frac{r_1^{n-2}(s)}{ u_1'( r_1(s))\sqrt{( u_1'( r_1(s))^2+2F(s)}}-
\frac{ r_2^{n-2}(s)}{u_2'(r_2(s))\sqrt{(u_2'(r_2(s)))^2+2F(s)}}\Bigr)\\
&\ge&r_2^{n-2}(s)|F(s)|\Bigl(
\frac{1}{| u_1'( r_1(s))|\sqrt{( u_1'( r_1(s))^2+2F(s)}}\\&&\quad -
\frac{ 1}{|u_2'(r_2(s))|\sqrt{(u_2'(r_2(s)))^2+2F(s)}}\Bigr) \\&\ge&0.
\end{eqnarray*}
\end{proof}

For the case when $f$ satisfies $(f_5)$-$(f_6)$ we use \cite[Proposition 4.1.2]{cghy}. Even though in this proposition we assumed $f$ superlinear, this assumption is not used in the proof, so we state it here without proof.

\begin{lema}\label{sep2mm}
Let $f$ satisfy $(f_1)$-$(f_2)$ and $(f_5)$-$(f_6)$. Then there exists $\delta\in(0,\delta_1]$ such that for all $\alpha_1,\ \alpha_2\in(\alpha^*-\delta,\alpha^*+\delta)$ with $\alpha_1<\alpha_2$ it holds that
$$r_1( s)> r_2( s)\quad\mbox{and}\quad \tilde W_1( s)<\tilde W_2( s),\quad\mbox{for all }s\in[-\beta,U_{bI}),$$
where $U_{bI}=\min\{b,U_I\}$.
\end{lema}

\begin{proof}[{\bf Proof of Proposition \ref{pmaxmin}}]
We prove this proposition in the case that $f$ satisfies $(f_1)$-$(f_3)$ and $(f_4)$, the proof when $f$ satisfies $(f_5)$-$(f_6)$ follows similarly by using Lemma \ref{sep2mm}. As in \cite{et, st}, we set
$$S_{12}(s)=\frac{r_1^{n-1}r_2'}{r_2^{n-1}r_1'}(s).$$
Then
\begin{equation}\label{sderiv}
  S'_{12}(s)=S_{12}(s)f(s)((r_2'(s))^2-(r_1'(s))^2).
\end{equation}
Let
$$U=\min\{-\beta, U_I\}.
$$

We will prove first that $m_{1}>m_{2}$ and that for all $s\in[m_{1},U)$ we have
\begin{equation} \label{conclusion} S_{12}(s)<1, \quad |r_1'(s)|>|r_2'(s)|, \quad r_1(s)>r_2(s).\end{equation}

If $U_I>-\beta$ then $U=-\beta$, and, from Lemma \ref{sep1}, and using that $F(-\beta)=0,$ we have that $S_{12}(U)\le 1$ and $r_1(U)>r_2(U)$. Thus, $|r_1'(U)|>|r_2'(U)|$. On the other hand, if $U=U_I$, we also have that $S_{12}(U)<1$ and $|r_1'(U)|>|r_2'(U)|$.

From \eqref{sderiv} we have that $S_{12}$(s) is increasing as long as $|r'_1(s)|>|r'_2(s)|$, for $s<U.$ If \eqref{conclusion} does not hold for all $s\in(\max\{m_{1},m_{2}\},U)$, then at the largest point $s_0$ where it fails, we must have that $|r'_1(s_0)|=|r'_2(s_0)|$ and $r_1(s_0)>r_2(s_0)$ implying that $S_{12}(s_0)>1$, a contradiction. Thus \eqref{conclusion} holds in $(\max\{m_{1},m_{2}\},U)$, and hence $m_{1}=\max\{m_{1},m_{2}\}$.

Next we prove that $P_1>P_2 $ in $[m_{1},U].$ From the definition of $P_1$ and $P_2$ we have
\begin{eqnarray*}\bigl(P_1-P_2\bigr)(U)&= & \Bigl(\frac{r_2^n}{(r'_2)^2}-\frac{r_1^n}{(r'_1)^2}\Bigr)(U)+2n\frac{F}{f}(U)\Bigl(\frac{r_1^{n-1}(U)}{|r_1'(U)|}-\frac{r_2^{n-1}(U)}{|r_2'(U)|}\Bigr)\\
&\ge&\Bigl(\frac{r_2^n}{(r'_2)^2}-\frac{r_1^n}{(r'_1)^2}\Bigr)(U)\\
&= & \Bigl(\frac{r_2^n}{(r'_2)^2}
\Bigl[ 1-S^2_{12}\frac{r^{n-2}_2}{r^{n-2}_1}\Bigr]\Bigr)(U)>0.\end{eqnarray*}

On the other hand, from  $(f_4)$ and \eqref{conclusion},
\begin{eqnarray*}
&&\bigl(P_1-P_2\bigr)'(s)=
(S_{12}(s)-1)\Bigl(n-2-2n\Bigl(\frac{F}{f}\Bigr)'(s)\Bigr)\frac{r_2^{n-1}}{r_2'}(s)
<0,
\end{eqnarray*}
implying that $P_1>P_2$ in $[m_{1},U].$
In particular, $P_1(m_{1})>P_2(m_{1}).$  Now, since $P'_2>0$, we have that $P_2(m_{1})>P_2(m_{2}),$ and thus $P_1(m_{1})>P_2(m_{2}),$ ending the proof of the proposition.

\end{proof}
The analogue of Lemma \ref{sep4ultimo} for the case $k=2$ can be found in \cite[Lemma 4.4.1]{cghy}, we state it below for the sake of completeness. Set
$$\bar W(s,\alpha)= \bar r(s,\alpha)\sqrt{(u'(\bar r(s,\alpha),\alpha))^2+2F(s)}, \quad s\in[m_1(\alpha),S(\alpha)),$$
where
 $$\bar S_j:= \sup\{s\in (m_j,  U_2(\alpha_j))\ :\ \quad (u_j'(\bar r_j(s)))^2+2F(s)>0\}.$$

\begin{lema}\label{sep5}
Assume that $f$ satisfies $(f_1)$-$(f_3)$, and let $\alpha^*\in\mathcal G_2$. Let $\alpha_1,\ \alpha_2\in(\alpha^*-\delta,\alpha^*+\delta)$ with $\alpha_1<\alpha_2$.
Assume that there exists $U\in[-\beta,0]$ such that
\begin{equation} \label{barcompa3}   r_1(U)\ge   r_2(U)\quad\mbox{and}\quad \bar W_1( U)< \bar W_2(U).\end{equation}
Then,    $$\bar S_1\le \bar S_2$$ and
$$\bar r_1(s)>\bar r_2(s), \quad W_1(s)<W_2(s),\quad\mbox{and \; $u'_1(\bar r_1(s))<u'_2(\bar r_2(s))$} \quad s\in(U, \bar S_1].$$
\end{lema}

We define

\begin{eqnarray}\label{sbderiv}
\begin{gathered}
\bar P(s,\alpha)=-2n\frac{F}{f}(s)\frac{\bar r^{n-1}(s,\alpha)}{\bar r'(s,\alpha)}
-\frac{\bar r^n(s,\alpha)}{(\bar r'(s,\alpha))^2}-2\bar r^n(s,\alpha)F(s),\\
\bar P'(s,\alpha)=\Bigl(n-2-2n\Bigl(\frac{F}{f}\Bigr)'(s)\Bigr)\frac{\bar r^{n-1}(s,\alpha)}{\bar r'(s,\alpha)},\\
\bar S_{12}(s)=\frac{\bar r_1^{n-1}\bar r_2'}{\bar r_2^{n-1}\bar r_1'}(s), \\
\bar S'_{12}(s)=\bar S_{12}(s)f(s)((\bar r_2'(s))^2-(\bar r_1'(s))^2).
\end{gathered}
\end{eqnarray}
\begin{prop}\label{sep5p}
Assume that $f$ satisfies $(f_1)$-$(f_3)$ and $(f_4)$, and let $\alpha^*\in\mathcal G_2$. Then there exists $\delta>0$ such that for $\alpha_1,\ \alpha_2\in(\alpha^*-\delta,\alpha^*+\delta)$ with $\alpha_1<\alpha_2$ it holds that:

if $\alpha_1\in \mathcal G_2\cup \mathcal N_2$,
then $\alpha_2\in \mathcal N_2$, 
\begin{equation}\label{sep55p}Z_2(\alpha_1)>Z_2(\alpha_2)\quad\mbox{and}\quad |u'_1(Z_2(\alpha_1))|<|u'_2(Z_2(\alpha_2))|,\end{equation}
and if $\alpha_2\in\mathcal G_2$, then $\alpha_1\in{\mathcal F_2}$.
\end{prop}

\begin{proof}
Let $m^*$ denote the minimum value of $u(\cdot,\alpha^*)$. Since $u'(r(m^*,\alpha^*),\alpha^*)=0$ and  $-2n\frac{F}{f}(m^*)>0$, by continuity we may choose $\delta\in(0,\delta_2)$ small enough so that
$$-2n\frac{F}{f}(m_{1})> \bar r_2(m_{1})u_2'(\bar r_2(m_{1})),$$
for all $\alpha_1,\ \alpha_2\in(\alpha^*-\delta,\alpha^*+\delta)$ and hence
\be\label{cnuevo}-2n\frac{F}{f}(m_{1})(\bar r_2(m_{1}))^{n-1}u_2'(\bar r_2(m_{1}))-(\bar r_2(m_{1}))^n(u_2'(\bar r_2(m_{1})))^2>0.\ee
On the other hand, from \eqref{ma1} in Proposition \ref{pmaxmin}, we have that $P_1(m_{1})>P_2(m_{2})$ and thus, using $m_{2}<m_{1}$ and the fact that $\bar P_2$ decreases, we find that
$$\bar P_1(m_{1})=P_1(m_{1})>P_2(m_{2})=\bar P_2(m_{2})>\bar P_2(m_{1}).$$ Therefore,
\ben
0&>&(\bar P_2-\bar P_1)(m_{1})\\
&=&-2n\frac{F}{f}(m_{1})(\bar r_2(m_{1}))^{n-1}u_2'(\bar r_2(m_{1}))-(\bar r_2(m_{1}))^n(u_2'(\bar r_2(m_{1})))^2\\
&&-2F(m_{1})(\bar r_2^n-\bar r_1^n)(m_{1})
\een
implying, by \eqref{cnuevo},
$$\bar r_1(m_{1})<\bar r_2(m_{1}).$$

We recall that from lemma \ref{delta1k}(ii), there exists an intersection point in $(T_1(\alpha^*),Z_2(\alpha^*))$. If $\bar r_{I}$ denotes the first of such points and if $\bar U_{I}= u_1(\bar r_{I})=u_2(\bar r_{I})$, then $\bar U_{I}\in(\overline{U}_1(\alpha^*),0]$.
Let us set
$$U=\max\{-\beta, \bar U_I\}.$$
We will show that $U$ satisfies \eqref{barcompa} in Lemma \ref{sep5}, that is,
\be\label{carmencita0} \bar r_1(U)\ge\bar r_2(U) , \quad \mbox{and}  \quad \bar W_1(U)<\bar W_2(U).\ee

We distinguish two cases:

\noindent Case 1.  $\bar U_I\in  [m_1,-\beta].$
We will first prove
\be\label{carmencita}\frac{\bar r_1^{n-1}}{\bar r_1'}(s)<\frac{\bar r_2^{n-1}}{\bar r_2'}(s)\quad\mbox{
and }\quad \bar P_1(s)>\bar P_2(s)\quad\mbox{ for all }s\in [m_{1},\bar U_I].\ee

Observe first that $\bar S_{12}(m_{1})=0$ and $\bar S_{12}(\bar U_I)<1$. If there exists a point $t\in(m_{1},\bar U_I)$ such that
$ \bar S'_{12}(t)=0,$ then $\bar r_1'(t)=\bar r_2'(t)$ and hence, from the definition of $\bar U_I$, $$ \bar S_{12}(t)= \frac{\bar r_1^{n-1}}{\bar r_2^{n-1}}(t)<1,$$ implying $\bar S_{12}(s)<1$ for $s\in [m_{1},\bar U_I]$.

On the other hand, from the second equation in \eqref{sbderiv}, using that $\bar S_{12}(s)<1$ and $(f_4)$, we obtain
\begin{eqnarray*}
&&\bigl(\bar P_1-\bar P_2\bigr)'(s)=
\Bigl ((\bar S_{12}-1)\Bigl(n-2-2n\Bigl(\frac{F}{f}\Bigr)'\Bigr)\frac{\bar r_2^{n-1}}{\bar r_2'}\Bigr)(s)
>0.
\end{eqnarray*}
Hence, for all $s\in (m_{1},\bar U_I),$ $\bar P_1(s)-\bar P_2(s)>\bar P_1(m_{1})-\bar P_2(m_{1})>0$

Next we will prove that
\begin{equation} \label{para} \bar r_1(s)>\bar r_2(s) , \quad \mbox{and}  \quad \frac{\bar r_1}{\bar r'_1}(s)<\frac{\bar r_2}{\bar r'_2}(s)\quad \mbox{for all $s\in(\bar U_{I}, -\beta]$. }
\end{equation}

From the definition of $\bar U_I$, $\displaystyle \frac{\bar r_1}{\bar r'_1}<\frac{\bar r_2}{\bar r'_2}$ at $\bar U_{I}$.
Assume by contradiction that \eqref{para} does not hold.  Then,  there exists a first point $t\in( \bar U_{I},-\beta)$ such that
$$\frac{\bar r_1}{\bar r'_1}(t)=\frac{\bar r_2}{\bar r'_2}(t)\quad \mbox{and } \bar r_1(s) >\bar r_2(s),\quad \mbox{for all $s\in (\bar U_{I}, t]$}, $$
implying
$$\bar{S}_{12}(t) = \Bigl(\frac{\bar r_1(t)}{\bar r_2(t)}\Bigr)^{n-2}= D >1.$$
From the definition of $\bar P_1$ and $\bar P_2$, we have that
$$
(\bar P_1-D\bar P_2)(t)= 2(D\bar r_2^{n}-\bar r_1^{n})F(t)=2\bar r_1^{n-2}(\bar r_2^2-\bar r_1^2)F(t)<0.
$$
On the other hand, from  \eqref{carmencita},   we have that $(\bar P_1-\bar P_2)(\bar U_{I})>0$.   Since $\bar P_2(m_{2})<0$ and  $\bar P_2$ decreases in $(m_{2},-\beta)$, we have that $\bar P_2(\bar U_I)<0$. Hence, as $D>1$,  we conclude that
$$(\bar P_1-D\bar P_2)(\bar U_{I})>0.$$
From the last equation in \eqref{sbderiv} we obtain  that $\bar S_{12}$ is increasing in $(\bar U_{I},t)$ implying that $\bar S_{12}(s)< D.$ Finally, using $(f_4)$ we deduce
$$
\bigl(\bar P_1-D\bar P_2\bigr)'(s)=
\Bigl ((\bar S_{12}-D)\Bigl(n-2-2n\Bigl(\frac{F}{f}\Bigr)'\Bigr)\frac{\bar r_2^{n-1}}{\bar r_2'}\Bigr)(s)>0$$
for all $s\in(\bar U_{I},t)$ and thus
$$
(\bar P_1-D\bar P_2)(t)>0,$$ a contradiction. Hence, \eqref{carmencita} follows, and, since $F(-\beta)=0$, also \eqref{carmencita0}.

\noindent Case 2.  $\bar U_I\in  [-\beta,0).$ In this case $U=\bar U_I$, and \eqref{carmencita0} trivially holds.

Hence, by Lemma \ref{sep5}, we have
$\bar S_1\le \bar S_2$,
$$ r_1(s)> r_2(s), \quad  \bar W_1(s)< \bar W_2(s),\quad\mbox{and \; $u'_1( r_1(s))<u'_2( r_2(s))$} \quad\mbox{for all } s\in(U,S_1].$$
If $\alpha_1\in\mathcal G_2\cup\mathcal N_2$, then $S_1=0$ implying $S_2=0$ and $\alpha_2\in\mathcal G_2\cup\mathcal N_2$. As
$Z_2(\alpha_1)=\bar r_1(0)>\bar r_2(0)=Z_2(\alpha_2)$ and $u_1'(Z_2(\alpha_1))<u_2'(Z_2(\alpha_2))$ we conclude that $\alpha_2\in\mathcal N_2$.

If $\alpha_2\in\mathcal G_2$, then $\bar S_2=0$. As $u_2'(Z_2(\alpha_2))=0$, we conclude that $\bar S_1<0$ implying $\alpha_1\in\mathcal F_2$.

\end{proof}

\begin{proof}[{\bf Proof of Theorem \ref{main1}} ]
Let $\alpha^*\in\mathcal G_2$ hence by Lemma \ref{sep5}, $(\alpha^*,\alpha^*+\delta)\subset
\mathcal N_2$. Let
$$\bar\alpha=\sup\{\alpha>\alpha^*\ :\ (\alpha^*,\alpha)\subset\mathcal N_2\}.$$
Assume $\bar\alpha<\infty$. Since $\mathcal P_2$ and $\mathcal N_2$ are open, we deduce that $\bar\alpha\in\mathcal G_2$. By Lemma \ref{sep5}, $(\bar\alpha-\delta,\bar\alpha)\subset{\mathcal F_2}$, a contradiction, and thus $(\alpha^*,\infty)\subset \mathcal N_2$. Hence, there exists at most one solution of \eqref{eq2} with exactly one sign change in $(0,\infty)$.

\end{proof}
\section{The Dirichlet problem}\label{b0}
We begin this section by noting that under assumptions $(f_1)$-$(f_3)$, there might be non uniqueness of the solutions to the Dirichlet problem \eqref{eq2d} in some balls, that is, for some values of $\rho>0$.
Indeed, assume that in addition to $(f_1)$-$(f_3)$, it holds that
\be\label{fsubnu}\liminf_{s\to\infty}\frac{F(s)}{s^2}=0.\ee
Then the results in \cite{fls} hold, and  in particular, there exists a ground state solution of \eqref{eq2}. Let $\alpha^*$ be the greatest initial value which gives rise to this solution. If the support of this solution is not compact, then for $\alpha>\alpha^*$ but close, it happens that $\alpha\in\mathcal N_1$ and $Z_1(\alpha)\to\infty$  as $\alpha\downarrow \alpha^*$. If the solution has compact support, then from Proposition \ref{varphi2} (ii), for $\alpha>\alpha^*$ but close enough, $Z_1(\alpha)<Z_1(\alpha^*)$.

On the other hand, by denoting by $r(\beta,\alpha)$ the first positive value of $r$ at which $u(r,\alpha)=\beta$. Since $F(u(r))>0$ for $r\in[0,r(\beta,\alpha)]$, we have that $|u'(r)|\le \sqrt{2F(\alpha)}$ for all $r\in[0,r(\beta,\alpha)]$. Hence, from the mean value theorem, there exists $\xi\in[0,r(\beta,\alpha)]$ such that
$$\frac{\alpha-\beta}{r(\beta,\alpha)}=|u'(\xi)|\le\sqrt{2F(\alpha)},$$
implying that
$$r(\beta,\alpha)\ge \frac{\alpha-\beta}{\sqrt{2F(\alpha)}},$$
and thus, from \eqref{fsubnu}, there exists a sequence $\alpha_i\to\infty$ as $i\to\infty$ such that $r(\beta,\alpha_i)\to\infty$ as $i\to\infty$.
From Proposition \ref{varphi1}, $r(\beta,\alpha)$ is increasing in $\alpha$, and hence
$\lim\limits_{\alpha\to\infty}r(\beta,\alpha)=\infty$. This in turn implies that
$$\lim_{\alpha\to\infty}Z_1(\alpha)=\infty,$$
and our claim follows.
\bigskip

Let now $f$ satisfy $(f_1')$-$(f_3')$ and $(f_4)$ (with $\beta=0$). We claim that there cannot exist bound state solutions to \eqref{eq2} with a finite number of zeros. Indeed, we first observe that from condition $(f_3')$, it easily follows that for any $s_0>0$ there exists a positive constant $C_0$ such that
\be\label{fsub}
f(s)\ge C_0s\quad\mbox{for all }s\in(0,s_0).
\ee
Now let $u$ be a solution to \eqref{eq2} with, say, $k$ zeros. Without loss of generality we may assume $0<u(r)<s_0$ for $r$ large, hence $u$ decreases for $r\ge r_0$, $r_0$ large, and thus
$$-r^{n-1}u'(r)\ge \int_{r/2}^rt^{n-1}f(u(t))dt\ge C_0\int_{r/2}^rt^{n-1}u(t)dt\ge Cr^nu(r)$$
for all $r\ge 2r_0$ and some positive constant $C$. Hence
$$-\frac{u'(r)}{u(r)}\ge Cr,\quad r\ge 2r_0,$$
and thus
\be\label{cotau}u(r)\le C_1e^{-r^2}\quad\mbox{for all $r\ge 2r_0$}\ee
and some positive constant $C_1$. Setting
$$P(r)=-2n\frac{F}{f}(u(r)) r^{n-1} u'(r)
-r^n(u'(r))^2-2 r^nF(u(r)),$$
we have that
$$P'(r)=\Bigl(n-2-2n\Bigl(\frac{F}{f}\Bigr)'(u(r))\Bigr)r^{n-1}|u'(r)|^2<0.$$
Since $P(0)=0$, and thanks to \eqref{cotau} and $(f_3')$, also $\lim\limits_{r\to\infty}P(r)=0$, we obtain a contradiction.

Hence, for $f$ satisfying $(f_1')$-$(f_3')$, we are led to  study the uniqueness of solutions with a prescribed number of zeros to the Dirichlet problem \eqref{eq2d}.

\noindent {\bf Proof of Theorem \ref{main2d0}.} It is based on the following facts:
\begin{enumerate}
\item Proposition \ref{varphi2} for the case $b=0$, that is, between two consecutive zeros of $u'$ there is at least one zero of $\varphi$.
\item The identity
$$(r^{n-1}(u'\varphi-\varphi'u))'=r^{n-1}\varphi(uf'(u)-f(u))$$
and condition $(f_3')$ say that there cannot be two zeros of $\varphi$ in $(Z_i(\alpha),Z_{i+1}(\alpha))$, and
\item We have that
$$\frac{d}{d\alpha}Z_i(\alpha)=-\frac{\varphi}{u'}>0\quad \mbox{for all $i$.}$$
\end{enumerate}
The details are left to the interested reader.

\end{document}